\documentclass[11pt]{article}

\usepackage{amsmath,amsthm}

\usepackage{amssymb,latexsym}

\usepackage{enumerate}

\newtheorem{tw}{Theorem}[subsection]
\newtheorem{lm}[tw]{Lemma}
\newtheorem{wn}[tw]{Corollary}
\newtheorem{stw}[tw]{Proposition}
\newenvironment{dow}{\it Proof.\rm}{\hfill $\Box$}

\theoremstyle{definition}
\newtheorem*{df}{Definition}
\newtheorem{uw}[tw]{Remark}
\newtheorem{prz}[tw]{Example}

\newcommand{\BN}{{\mathbb N}}
\newcommand{\BR}{{\mathbb R}}
\newcommand{\BX}{{\mathbb X}}

\newcommand{\FF}{{\mathcal{F}}}

\newcommand{\HH}{{\mathcal{H}}}

\newcommand{\MM}{{\mathcal{M}}}

\newcommand{\EE}{{\mathcal{E}}}

\newcommand{\intot}{{\int_{0}^{t}}}

\newcommand{\BRD}{{\mathbb{R}^{d}}}

\newcommand{\nsubsection}{\setcounter{equation}{0}\subsection}

\begin{document}

\title {Semilinear elliptic systems with measure data}
\author {Tomasz Klimsiak \smallskip\\
{\small Faculty of Mathematics and Computer Science,
Nicolaus Copernicus University} \\
{\small  Chopina 12/18, 87--100 Toru\'n, Poland}\\
{\small e-mail: tomas@mat.uni.torun.pl}}
\date{}
\maketitle
\begin{abstract}
We study the Dirichlet problem for systems of the form $-\Delta
u^k=f^k(x,u)+\mu^k$, $x\in\Omega$, $k=1,\dots,n$,  where
$\Omega\subset\BR^d$ is an open (possibly nonregular) bounded set,
$\mu^1,\dots,\mu^n$ are bounded diffuse measures on $\Omega$,
$f=(f^1,\dots,f^n)$ satisfies some mild integrability condition
and the so-called angle condition. Using the methods of
probabilistic Dirichlet forms theory we show that the system has a
unique solution in the generalized Sobolev space
$\dot{H}^{1}_{loc}(\Omega)$ of functions having fine gradient. We
provide also a stochastic representation of the solution.
\end{abstract}


\footnotetext{{\em Key words or phrases:} Semilinear elliptic
systems, Laplacian, measure data, Dirichlet form.}

\nsubsection{Introduction}

Let $\Omega\subset\BRD$, $d\ge 2$, be an open  bounded set.
In the present paper we study the
existence and uniqueness of solutions of systems of the form
\begin{equation}
\label{eq1.1} \left\{
\begin{array}{l} -\frac12\Delta u^{k}
=f^{k}(x,u)+\mu^{k}\quad \mbox{in }\Omega,\quad k=1,\dots,n,\medskip\\
u^{k}=0\quad \mbox{on }\partial\Omega,\quad k=1,\dots,n,
\end{array}
\right.
\end{equation}
where
$f=(f^1,\dots,f^n):\Omega\times\mathbb{R}^{n}\rightarrow\mathbb{R}^n$
is a Carath\'eodory function and $\mu^1,\dots,\mu^n$ belong to the
space $\MM_{0,b}$ of bounded diffuse measures on $\Omega$ (see
Section 2).

Let $qL^{1}_{loc}(\Omega)$ denote the space of locally
quasi-integrable functions (see Section 2). In the scalar case
($n=1$) it is known that if
\begin{equation}
\label{eq.i2} x\mapsto \sup_{|u|\le r}|f(x,u)|\in
qL_{loc}^{1}(\Omega) \mbox{ for every } r \ge 0,
\end{equation}

\begin{equation}
\label{eq.i3} y\mapsto f(x,y) \mbox{ is continuous on } \mathbb{R}
\mbox{ for a.e. }x\in\Omega
\end{equation}
and
\begin{equation}
\label{eq.i4} f(x,u)\cdot u\le 0\mbox{ for a.e. } x\in\Omega
\mbox{ and every }u\in\mathbb{R},
\end{equation}
then there exists a solution of (\ref{eq1.1}) (see \cite{OP}; see
also \cite{Betal.} for equations with general Leray-Lions type
operators). One of the crucial ingredient in the proof of the
existence result for (\ref{eq1.1}) is the following Stampacchia
estimate
\begin{align}
\label{eq0.s}
\|f(\cdot,u)\|_{L^{1}(\Omega;m)}\le\|\mu\|_{TV}
\end{align}
(see \cite{Stampacchia}), which holds true under (\ref{eq.i4}) for
every solution of (\ref{eq1.1}). An attempt to generalize the
existence result to $n=2$ has been made in \cite{OP}. It is proved
there that if $\Omega$ is smooth, $f$ does not depend on $x$,  is
continuous on $\mathbb{R}^{2}$ and monotone componentwise, i.e.
$f^{1}(\cdot,v),\, f^{2}(u,\cdot)$ are nonincreasing and
$f^{1}(0,v)=f^{2}(u,0)=0$ for every $(u,v)\in\mathbb{R}^{2}$, then
there exists a unique solution of (\ref{eq1.1}). In \cite{OP}, as
in the scalar case, the key step in proving the existence of
solutions is Stampacchia's estimate, which  is derived by using
the componentwise character of monotonicity of $f$ and by
introducing the important notion of quasi-integrability of
functions. Note also that in \cite[Remark 7.1]{OP} the authors
raise the question of existence of solutions to (\ref{eq1.1}) for
$f$ satisfying weaker than monotonicity sign condition with
respect to each coordinate, i.e. for $f$ such that
$f^{1}(\cdot,u)\cdot v\le 0$, $f^{2}(u,\cdot)\cdot u\le 0$ for
every $(u,v)\in\mathbb{R}^{2}$. We answer positively the question
raised in \cite{OP}. Actually, using quite different than in
\cite{OP} methods of proof we show existence and uniqueness
results for more general systems.

In the present paper we assume that $\mu$, $f$ satisfy the
following assumptions analogous to  assumptions
(\ref{eq.i2})--(\ref{eq.i4}) considered in the theory of scalar
equations:
\begin{enumerate}
\item[(A1)] $\mu\in\MM_{0,b}$,
\item[(A2)]For every  $r\ge 0$,
$x\mapsto\sup_{|y|\le r}|f(x,y)|\in q{L}^{1}_{loc}(\Omega)$,
\item[(A3)]For a.e. $x\in\Omega$, $y\rightarrow f(x,y)$ is continuous,
\item[(A4)]$\langle f(x,y),y\rangle\le 0$, $y\in\mathbb{R}^{n}$
and a.e. $x\in\Omega$ (Here $\langle\cdot,\cdot\rangle$ stands for
the usual scalar product in $\BR^n$).
\end{enumerate}
Following \cite{Landes} we will call (A4) the angle condition. In
\cite{Landes} more general than (\ref{eq1.1}) elliptic systems
with the perturbed Leray-Lions type operator are considered.
As a matter of fact the assumptions in \cite{Landes} when adjusted
to our problem say that the perturbation satisfies some strong
growth conditions and stronger then (A4) condition
\begin{enumerate}
\item[(A5)]There exists $\alpha>0$ such that for every
$y\in\mathbb{R}^{n}$ and a.e. $x\in\BRD$,
\[
\langle f(x,y),y\rangle \le -\alpha|y|^2,
\]
\end{enumerate}
which we will call the uniform angle condition.

In the paper we show that if the right-hand side of (\ref{eq1.1})
satisfies (A5) then Stampacchia's estimate (\ref{eq0.s}) holds
true for any solution of (\ref{eq1.1}), which immediately implies
that any solution of (\ref{eq1.1}) belongs to the Sobolev space
$W^{1,q}_{0}(\Omega)$ with $0<q<\frac{d}{d-1}$. Under (A4) no
analogue of Stampacchia's estimate appears to be available.
Consequently, it seems that in general $f(\cdot,u)\notin
L^{1}_{loc}(\Omega)$. Therefore the first problem we have to
address is to describe the regularity space for $u$ and then
formulate suitable definition of a solution of (\ref{eq1.1}). We
propose two equivalent definitions: the probabilistic and analytic
one.

Let $\BX=(X,P_x)$ denote the Wiener process killed upon leaving
$\Omega$ and let $\zeta$ denote its life-time. In the
probabilistic definition by a solution we mean a quasi-continuous
in the restricted sense function $u:\Omega\rightarrow\BR$ such
that the following stochastic equation
\begin{align}
\label{eq1.6}
u(X_{t})&=u(X_{\tau})+\int_{t}^{\tau}f(X_{\theta},u(X_{\theta}))\,d\theta
\nonumber\\
&\quad+\int_{t}^{\tau}dA^{\mu}_{\theta}
-\int_{t}^{\tau}dM_{\theta}, \quad 0\le t\le\tau,\quad
P_x\mbox{-a.s.}
\end{align}
is satisfied for quasi-every (with respect to the Newtonian
capacity) $x\in\Omega$. Here $\tau$ is an arbitrary stopping time
such that $0\le\tau<\zeta$, $M$ is some local martingale additive
functional of $\BX$ (as a matter of fact $M$ corresponds to the
gradient of $u$) and $A^{\mu}$ is a positive continuous additive
functional of $\BX$ associated with the measure $\mu$ via Revuz
duality (see Section \ref{sec2}).

In the analytic definition (see Section \ref{sec4}), a
quasi-continuous function $u:\Omega\rightarrow\BR$ is a solution
of (\ref{eq1.1}) if
\begin{equation}
\label{eq1.7} \EE(u^{k}, v)=(f^{k}(\cdot,u),v)_{{L}^{2}(\Omega;m)}
+\langle \mu,v\rangle, \quad v\in H^{1}_{0}(G_{l})
\end{equation}
for every $k=1,\dots,n$ and $l\ge0$, where
\begin{equation}
\label{DF} \EE(u,v)=\frac12\int_{\Omega}\nabla u(x)\nabla
v(x)\,m(dx), \quad u,v\in D[\EE]=H^{1}_{0}(\Omega)
\end{equation}
is the Dirichlet form associated with the operator
$(\frac12\Delta, H^{1}_{0}(\Omega))$, $\{G_{l},\,l\ge 1\}$ is a suitable
family of finely open sets depending on $u$ such $\bigcup_{l\ge 1} G_{l}=\Omega$
q.e. and (\ref{eq1.7}) makes sense.

Since we are looking for solutions of (\ref{eq1.1}) in Sobolev
type spaces, our minimal regularity requirement for them is
quasi-continuity or, equivalently, continuity in the fine topology
(see \cite{Fukushima}). Quasi-continuity provides some information
on local, in terms of fine topology, regularity of functions and
allows one to control their behavior on finely open (closed) sets
of the form $\{u<t\}$ ($\{u\le t\}$), $t\ge 0$.
It is therefore natural to try to derive a priori estimates for
solutions of (\ref{eq1.1}) in Sobolev spaces on finely open sets
to make sense of the analytic definition and then prove existence
of solutions of (\ref{eq1.1}). In the present paper we prefer,
however, a stochastic approach to the problem. The main reason for
adopting the stochastic approach is that (\ref{eq1.6}) is simpler
to investigate than (\ref{eq1.7}), because (\ref{eq1.7}) is in
fact a family of variational equations on finely open domains
$G_{k}$ which depend on the solution. Equations of the form
(\ref{eq1.7}) were considered for example in \cite{Fuglede2,KM}.
It seems that direct analysis of equations (\ref{eq1.7}) would
generate many technical difficulties in using the fine topology,
while the stochastic approach avoids them, because in the latter
approach the fine topology is hidden in a very convenient way in
probabilistic notions of the Dirichlet forms theory we are using
in our proofs.

We prove that under (A1)--(A4) there exists a probabilistic
solution $u$ of (\ref{eq1.1}), $f(\cdot,u)\in qL^{1}_{loc}(\Omega)$ and
$u$ belongs to the generalized Sobolev space
$\dot{H}^{1}_{loc}(\Omega)$ of functions having fine gradient (see
\cite{KM}). The space $\dot{H}^{1}_{loc}(\Omega)$ is wider than
the space $\mathcal{T}^{1,2}$ of  Borel functions whose
truncations on every level belong to $H^{1}_{0}(\Omega)$, which
was introduced in \cite{Betal.} to cope with elliptic equations
with $L^{1}$ data. The solution $u$ satisfies (\ref{eq1.7}),
because we show that in general, if
$u\in\dot{H}^{1}_{loc}(\Omega)$ then $u$ is a probabilistic
solution of (\ref{eq1.1}) iff it is a solution of (\ref{eq1.1}) in
the analytic sense. We also prove that if we replace (A4) by
\begin{enumerate}
\item [(A4$'$)]$\langle f(x,y)-f(x,y'), y-y'\rangle\le 0
\mbox{ for a.e. }x\in\Omega\mbox{ and every }
y,y'\in\mathbb{R}^{n}$,
\end{enumerate}
then the solution of (\ref{eq1.1}) is unique. Let us note that
besides proving our existence result under the general angle
condition (A4), in contrast to \cite{OP} we allow $f$ to depend on
$x$, the dimension of the system is arbitrary and we impose no
assumption on the regularity of $\Omega$.

Moreover, we show that if $u$ is a solution of (\ref{eq1.1})  and
$f(\cdot,u)\in L^{1}(\Omega;m)$ then $u\in W^{1,q}_{0}(\Omega)$ with
$1\le q<d/(d-1)$ and $u$ coincides
with the distributional (renormalized, in the sense of duality)
solution. Finally, we show that $f(\cdot,u)\in L^{1}(\Omega;m)$ if
\begin{enumerate}
\item[(A4$''$)] For every $k\in\{1,\dots,n\}$,
$y_{i}\in \mathbb{R}$, $i\in\{1,\dots,n\}\setminus\{k\}$ and a.e.
$x\in\Omega$,
\[
f^{k}(x,y_{1},\dots,y_{k-1},z,y_{k+1},\dots,y_{n})\cdot z\le 0
\mbox{ for every }z\in\mathbb{R},
\]

\end{enumerate}
i.e. if $f$ satisfies the sign condition which respect to each
coordinate.


\nsubsection{Preliminary results} \label{sec2}

In the whole paper we adopt the convention that for given class of
real functions (measures) $F$ and $\mathbb{R}^{d}$-valued function (measure) $f$ we
write $f\in F$ if each component $f^{k}$, $k=1,\dots,d$, of $f$
belongs to $F$.

Let $\Omega$ be a nonempty bounded open subset in $\BRD$, $d\ge2$.
By $m$ we denote the Lebesgue measure on $\BRD$.
Let $(\EE,D[\EE])$ be the Dirichlet form defined by (\ref{DF}) and
let $\mathbb{X}=(\{X_{t},\,t\ge 0\},\{P_{x},
x\in\Omega\},\{\FF_{t},\,t\ge 0\},\zeta)$ be a diffusion process
associated with $(\EE,D[\EE])$, i.e. for  every $t\ge 0$,
\begin{equation}
\label{SR} (p_{t}f)(x)=E_{x}f(X_{t})
\end{equation}
for $m$-a.e. $x\in\Omega$, where $\{p_{t},t\ge 0\}$ is the
semigroup generated by $(\EE, D[\EE])$ and $\zeta=\inf\{t\ge 0,
X_{t}=\Delta\}$, where $\Delta$ is a one-point compactification of
the space $\Omega$ and $E_{x}$ denote the expectation with respect
to $P_{x}$ (see \cite{Fukushima}). We also admit the convention
that $u(\Delta)=0$. It is well known that $\mathbb{X}$ is the Brownian motion
killed upon leaving $\Omega$ so we sometimes  use the letter $B_{t}$ instead of $X_{t}$.
From (\ref{SR}) it follows
that $\{p_{t},\, t\ge 0\}$ is a semigroup of contractions.

By $\{R_{\alpha},\alpha>0\}$ we denote the resolvent generated by
$\{p_{t},\, t\ge 0\}$. Since $\EE$ is transient (see \cite[Example
1.5.1]{Fukushima}), $R_{0}$ is also well defined. Write $R=R_{0}$.
From (\ref{SR}) it follows that
\[
(Rf)(x)=E_{x}\int_{0}^{\zeta} f(X_{t})\,dt
\]
for $m$-a.e. $x\in\Omega$.

By $\mbox{cap}$ we denote the capacity associated with $\EE$, i.e.
$\mbox{cap}: 2^{\Omega}\rightarrow \mathbb{R}^{+}\cup \{+\infty\}$
is a subadditive set functions defined as
\[
\mbox{cap}(A)=\inf\{\EE(u,u); u\in
H^{1}_{0}(\Omega),\,u\ge\mathbf{1}_{A},\,\, m\mbox{-a.e.}\}
\]
for open set $A\subset\Omega$ (with the convention that
$\inf\emptyset=+\infty$), and for arbitrary $A\subset\Omega$,
$\mbox{cap}(A)=\inf\{\mbox{cap$(B):A\subset B,\,B$ is an open
subset of $\Omega$}\}$. We say that some property $P(x)$ holds
q.e. if $\mbox{cap}(\{x; P(x)\mbox{ is not true}\})=0$.

A Borel measurable set $N\subset \Omega$ is called properly
exceptional if for q.e. $x\in\Omega$,
\begin{equation}
\label{exceptional}
P_{x}(\exists t>0; X_{t}\in N)=0.
\end{equation}
It is known (see \cite[Theorem 4.1.1 and page 140]{Fukushima})
that every properly exceptional set is of capacity zero and if
$N\subset \Omega$ is of capacity zero then there exists a Borel
properly exceptional set $B$ such that $N\subset B$.

A nonnegative Borel  measure $\mu$ on $\Omega$ is called smooth if
it charges no set of zero capacity and there exists an ascending
sequence $\{F_{n}\}$ of closed subsets of $\Omega$ such that
$\mu(F_{n})<\infty$ for $n\ge 1$ and for every compact set
$K\subset \Omega$,
\begin{equation}
\label{van} \mbox{cap}(K\setminus F_{n})\rightarrow 0.
\end{equation}
By $S$ we denote the set of all smooth measures on $\Omega$. By
$\MM^{+}_{0,b}$ we denote the space of finite smooth measures,
$\MM_{0,b}=\MM^{+}_{0,b}-\MM^{+}_{0,b}$.  Elements of $\MM_{0,b}$
are called diffuse or soft measures (see \cite{DMOP,DPP}).

Let $\mathcal{B}(\Omega)$ ($\mathcal{B}^+(\Omega)$) denote the set
of all real (nonegative) Borel measurable functions on $\Omega$.
For $A\subset \Omega$ we write $A\in\mathcal{B}(\Omega)$ if
$\mathbf{1}_{A}\in\mathcal{B}(\Omega)$. It is known (see
\cite[Section 5.1]{Fukushima}) that for every $\mu\in S$ there
exists a unique positive continuous additive functional $A^{\mu}$
(PCAF for short) such that for every
$f,h\in\mathcal{B}^{+}(\Omega)$,
\begin{equation}
\label{RD1}
E_{h\cdot m}\int_{0}^{t}f(X_{\theta})\, dA^{\mu}_{\theta} =\intot
\langle f\cdot\mu,p_{\theta} h\rangle \,d\theta,\quad t\ge 0,
\end{equation}
where $P_{\nu}(B)=\int_{\Omega} P_{x}(B)\,\nu(dx)$ for any
$B\in\mathcal{B}(\Omega)$ and nonnegative Borel measure $\nu$,
$(h\cdot\nu)(B)=\int_{B} h(x)\,\nu(dx)$ for any
$B\in\mathcal{B}(\Omega)$ and $h\in\mathcal{B}^{+}(\Omega)$, and
\[
\langle h,\nu\rangle =\int_{\Omega} h(x)\,\nu(dx).
\]
On the other hand, for every PCAF $A$ of $\mathbb{X}$ there exists
a unique smooth measure $\mu$ such that (\ref{RD1}) holds with
$A^{\mu}$ replaced by $A$. The measure $\mu$ is called the Revuz
measure associated with PCAF $A$.

Using (\ref{RD1}) one can extend the resolvent $R$ to $S$ by
putting
\begin{equation}
\label{RD2} R\mu (x)=E_{x}\int_{0}^{\zeta}dA_t^{\mu}.
\end{equation}
Similarly one can extend $R_{\alpha},\alpha>0$. The right-hand
side of (\ref{RD2}) is finite q.e. for every finite measure
$\mu\in S$. This follows from (\ref{TV}) and Lemma 4.1 and
Proposition 5.13 in \cite{KR}. Moreover, by Theorem 2.2.2,  Lemma
2.2.11 and Lemma 5.1.3 in \cite{Fukushima}, for every $\mu\in
H^{-1}(\Omega)$ and $v\in H^{1}_{0}(\Omega)$,
\begin{equation}
\label{dualitymeasure}
\EE(R\mu, \tilde{v})=\langle \mu, \tilde{v}\rangle,
\end{equation}
where $\tilde{v}$ is a quasi-continuous $m$-version of $v$.

By $\mathcal{C}(\Omega)$ we denote the space of all
quasi-continuous functions on $\Omega$. Let us recall that
$u:\Omega\rightarrow \mathbb{R}$ is quasi-continuous if for every
$\varepsilon>0$ there exists an open set $G_{\varepsilon}\subset
\Omega$ such that $\mbox{cap}(G_{\varepsilon})<\varepsilon$ and
$u_{|\Omega\setminus G_{\varepsilon}}$ is continuous. It is known
that $u$ is quasi-continuous iff the process $t\rightarrow
u(X_{t})$ is continuous on $[0,\zeta)$, $P_{x}$-a.s. for q.e.
$x\in\Omega$ (see \cite[Section 4.2]{Fukushima}).

$\mathcal{C}_{0}(\Omega)$ is the space of quasi-continuous
functions on $\Omega$ in the restricted sense. Let us note that a
Borel measurable function $u$ on $\Omega$ when considered as the
function on $\Omega\cup\{\Delta\}$ (with the convention that
$u(\Delta)=0$) is continuous on $\Omega\cup \{\Delta\}$ if $u\in
C_{0}(\Omega)$, where $C_{0}(\Omega)$ is the closure in
$C(\overline{\Omega})$ of the space $C_{c}(\Omega)$ of all
continuous functions with compact support in $\Omega$. A Borel
function $u$ on $\Omega$ is called quasi-continuous in the
restricted sense if for every $\varepsilon>0$ there exists an open
set $G_{\varepsilon}\subset\Omega$ such that
$\mbox{cap}(G_{\varepsilon})<\varepsilon$ and
$u_{|(\Omega\cup\{\Delta\})\setminus G_{\varepsilon}}$ is
continuous.

The following lemma shows that if $u\in\mathcal{C}_{0}(\Omega)$
then $u(x)$ tends to zero if $x$ tends to the boundary of $\Omega$
along the trajectories of the process $X$.
\begin{lm}
Assume that $u\in \mathcal{C}_{0}(\Omega)$. Then for q.e.
$x\in\Omega$,
\[
u(X_{t})\rightarrow 0,\quad t\rightarrow \zeta^{-},\quad
P_{x}\mbox{\rm-a.s.}
\]
\end{lm}
\begin{dow}
Since $u\in\mathcal{C}_{0}(\Omega)$ there exists a sequence
$\{U_{n}\}$ of open subsets of $\Omega$ such that
$\mbox{cap}(U_{n})\rightarrow 0$ and
$u_{|\Delta\cup\Omega\setminus U_{n}}$ is continuous. For every
$T>0$,
\begin{align*}
P_{x}(u(X_{t})\nrightarrow 0,\,t\rightarrow \zeta^{-}) &\le
P_{x}(T\le\zeta)+P_{x}(u(X_{t})\nrightarrow 0,\,t\rightarrow
\zeta^{-}, \zeta\le T)\\
& \le P_{x}(T\le\zeta)+ P_{x}(\exists_{t\le T}\, X_{t}\in U_{n})\\
&\le P_{x}(T\le\zeta)+P_{x}(\sigma_{U_{n}}\le T)\\
& \le P_{x}(T\le\zeta)+e^TE_{x}e^{-\sigma_{U_{n}}},
\end{align*}
where
\[
\sigma_{U_{n}}=\inf\{t>0;\,\, X_{t}\in U_{n}\}.
\]
Letting $n\rightarrow +\infty$ and applying \cite[Theorem
4.2.1]{Fukushima} we conclude that for q.e. $x\in\Omega$,
\[
P_{x}(u(X_{t})\nrightarrow 0,\,t\rightarrow \zeta^{-})\le
P_{x}(T\le\zeta).
\]
Since $E_{x}\zeta<+\infty$ (see (\ref{eq.esp})) for every
$x\in\Omega$, the desired result follows.
\end{dow}


\begin{uw}
\label{van.b} Let us note that by \cite[Theorem 2.1.3]{Fukushima},
if $u\in H_{0}^{1}(\Omega)$ then $u\in\mathcal{C}_{0}(\Omega)$.
\end{uw}

By $q{L}^{1}(\Omega)$ (resp. $q{L}^{1}_{loc}(\Omega)$) we denote
the class of all  Borel measurable  functions $f:\Omega\rightarrow
\mathbb{R}$ such that for q.e. $x\in\Omega$,
\begin{equation}
\label{qL}
P_{x}(\int_{0}^{T\wedge\zeta}|f(X_t)|\,dt<+\infty,\,T\ge0)=1
\end{equation}
\[
\left(\mbox{resp. } P_{x}(\int_{0}^{T}|f(X_t)|\,dt<+\infty,\,0\le
T <\zeta)=1 \right).
\]
Elements of $q{L}^{1}(\Omega)$ ($q{L}^{1}_{loc}(\Omega)$) will be
called quasi-integrable functions (locally quasi-integrable
functions).

We say that a Borel measurable function $f:\Omega\rightarrow
\mathbb{R}$ is locally quasi-integrable in the analytic sense if
for every compact $K\subset U$ and every $\varepsilon>0$ there
exists an open set $U_{\varepsilon}\subset \Omega$ such that
$\mbox{cap}(U_{\varepsilon})<\varepsilon$ and $f_{|K\setminus
U_{\varepsilon}}\in {L}^{1}(K\setminus U_{\varepsilon})$. We say
that a Borel measurable function $f:\Omega\rightarrow \mathbb{R}$
is quasi-integrable in the analytic sense if in the above
definition one can replace $K$ by $\Omega$.

\begin{uw}
From \cite{KR} it follows (see Remark 4.4) that $f$ is
quasi-integrable in the analytic sense iff $f\in
q{L}^{1}(\Omega)$ (this is true for bounded domains).
Moreover, if $f$ is locally quasi-integrable in
the analytic sense then $f\in q{L}^{1}_{loc}(\Omega)$. The reverse
implication is also true. Indeed, let $f\in
q{L}^{1}_{loc}(\Omega)$ and $f\ge 0$. Then by the very definition
of the space $q{L}^{1}_{loc}(\Omega)$, $A_{t}\equiv
\int_{0}^{t}f(X_{\theta})\,d\theta$ is a PCAF of $\mathbb{X}$ and
its associated Revuz measure is $f\cdot m$ (see \cite[Section
5.1]{Fukushima}). Since the associated Revuz measure is smooth,
there exists an ascending sequence $\{F_{n}\}$ of closed subsets
of $\Omega$ such that $f_{|F_{n}}\in{L}^{1}(F_{n};m)$ for every
$n\ge 1$ and (\ref{van}) holds for every compact set
$K\subset\Omega$. From this one can easily deduce that $f$ is
locally quasi-integrable in the analytic sense.
\end{uw}

The notion of quasi-integrability in the analytic sense was
introduced in \cite{OP}. In \cite{OP} the authors do not
distinguish between local quasi-integrability and
quasi-integrability, and quasi-integrabiity in the sense of
\cite{OP} coincides with local quasi-integrability in the
analytical sense defined in the present paper.

By $\mathcal{T}$ we denote the set of all stopping times with
respect to the filtration $\{\FF_{t},\,t\ge 0\}$ (see (\ref{SR})).
Let us recall (see \cite[Section 5.2]{Fukushima} that $M$ is
called a martingale additive functional (MAF) of $\BX$ if $M$ is a
finite continuous additive functional of $\mathbb{X}$ such that
for every $t>0$, $E_{x}M^{2}_{t}<\infty$ and $E_{x}M_{t}=0$ for
q.e. $x\in\Omega$. By $\MM$ we denote the space of all MAFs of
$\BX$. By $\MM_{loc}$ we denote the set of all local additive
functionals of $\BX$ (see \cite[page 226]{Fukushima}) for which
there exists an ascending sequence $\{G_{n},\, n\ge 1\}$ of finely
open subsets of $\Omega$ such that $\bigcup_{n\ge 1} G_{n}=\Omega$
q.e.,  a sequence $\{M^{n}\}\subset\MM$ and $N\subset\Omega$ such
that cap$(N)$=0 and for every $n\ge1$ and $x\in\Omega\setminus N$,
\[
M_{t}=M^{n}_{t},\quad t<\tau_{G_{n}},\quad P_{x}\mbox{-a.s.}
\]
Finally by $\mathcal{M}^{2}$ we denote the space of all $M\in\mathcal{M}$
such that $\sup_{t\ge 0}EM^{2}_{t}<\infty$ for q.e. $x\in\Omega$.

From now on we admit the following notation
\[ f(u)(x)=f(x,u(x)), \quad x\in\Omega\]
for every measurable function $u:\Omega\rightarrow \mathbb{R}^{n}$.

Following \cite{Kl3,KR}  let us consider the class (FD) consisting
of all functions  $u\in\mathcal{B}(\Omega)$ with the property that
the process $t\rightarrow u(X_{t})$ is of Doob's class (D) under
the measure $P_{x}$ for q.e. $x\in\Omega$, i.e. for q.e.
$x\in\Omega$ the family $\{u(X_{\tau}),\tau\in\mathcal{T}\}$ is
uniformly integrable under $P_{x}$.

\begin{df}
We say that  $u:\Omega\rightarrow \mathbb{R}^n$ is a probabilistic
solution of (\ref{eq1.1}) if
\begin{enumerate}
\item [(a)] $u$ is of class (FD),
\item [(b)]$u\in\mathcal{C}_{0}(\Omega)$,
\item [(c)]$x\mapsto f(u)(x)$ belongs to $q{L}^{1}_{loc}(\Omega)$,
\item [(d)]There exists  $M\in\MM_{loc}$ such
that for every stopping time $0\le\tau <\zeta$,
\begin{align}
\label{eq.df} u(X_{t})&=u(X_{\tau})+\int_{t}^{\tau}f(u)(X_{\theta})\,d\theta
+\int_{t}^{\tau}dA^{\mu}_{\theta}\nonumber \\
&\quad-\int_{t}^{\tau} dM_{\theta},\quad 0\le t\le
\tau,\quad P_x\mbox{-a.s.}
\end{align}
for q.e. $x\in\Omega$.
\end{enumerate}
\end{df}

In the sequel we admit the convention that $\int_{a}^{b}=0$ if
$a\ge b$.

\begin{uw}
\label{uw2.4} Under (b), if moreover $f(u)\in q{L}^{1}(\Omega)$,
(\ref{eq.df}) is satisfied iff for every $T\ge 0$,
\begin{align} \label{eq.df.equiv}
u(X_{t})&=u(X_{T\wedge\zeta})
+\int_{t}^{T\wedge\zeta}f(u)(X_{\theta})\,d\theta
+\int_{t}^{T\wedge\zeta}
dA^{\mu}_{\theta}\nonumber \\
&\quad-\int_{t}^{T\wedge\zeta} dM_{\theta},\quad t\ge 0,\quad P_x\mbox{-a.s.}
\end{align}
for q.e. $x\in\Omega$. Indeed, by (\ref{RD2}),
$E_{x}A_{\zeta}<\infty$ for q.e. $x\in\Omega$. Since every
stopping time with respect to Brownian filtration is predictable,
there exists a sequence $\{\tau_{n}\}\subset\mathcal{T}$ such that
$0\le\tau_{n}<T\wedge\zeta$ and $\tau_{n}\nearrow T\wedge\zeta$.
Taking $\tau_{n}$ in place of $\tau$ in (\ref{eq.df}) and letting
$n\rightarrow +\infty$ we get (\ref{eq.df.equiv}) (the integral
involving $f(u)$ is well defined since $f(u)\in
q{L}^{1}_{loc}(\Omega)$). On the other hand, if
(\ref{eq.df.equiv}) is satisfied, then replacing $T$ by an
arbitrary stopping time $\tau$ such that $0\le\tau<\zeta$ we get
(\ref{eq.df}).
\end{uw}

\begin{uw}
\label{uw.vanishing} If $f$ satisfies (A4) and $u$ is a solution
of (\ref{eq1.1}) then $u$ vanishes on the boundary of $\Omega$ in
the sense of Sobolev spaces. Indeed,  by the It\^o-Tanaka formula
(see \cite{BDHPS}) and (A4), for any $\tau\in\mathcal{T}$ such
that $0\le\tau<\zeta$,
\begin{align}
\label{eq2.tf}
\nonumber
E_{x}|u(X_{t})|&\le E_{x}|u(X_{T\wedge\zeta})|
+E_{x}\int_{t}^{T\wedge\zeta}
\langle f(u)(X_{\theta}),\hat{u}(X_{\theta})\rangle\,d\theta
+E_{x}\int_{t}^{T\wedge\zeta}
\langle\hat{u}(X_{\theta}),dA^{\mu}_{\theta}\rangle\\
&\le E_{x}\int_{0}^{\zeta}d|A^{\mu}|_t+E_{x}|u(X_{T\wedge\zeta})|
\end{align}
for q.e. $x\in\Omega$, where
\[
\hat{y}=\frac{y}{|y|}\mathbf{1}_{\{y\neq 0\}},\quad y\in\BR^d.
\]
Let $\{\tau_{k}\}$ be a sequence of
stopping times such that $0\le\tau_{k}<\zeta$, $k\ge 1$, and
$\tau_{k}\rightarrow \zeta$. Such a sequence exists since every
stopping time with respect to Brownian filtration is predictable
(see \cite[Theorem 4, Chapter 3]{Protter}).
It is clear that $u(X_{\tau_{k}})\rightarrow 0$ as $k\rightarrow
+\infty$, $P_{x}$-a.s.  for q.e. $x\in\Omega$. This when combined
with the fact that $u$ is of class (FD) implies that
\begin{equation}
\label{eq2.8} |u(x)|\le E_{x} \int_{0}^{\zeta}d|A^{\mu}|_t\equiv
v(x)
\end{equation}
for q.e. $x\in\Omega$. By \cite{KR}, $T_{k}(v)\in
H^{1}_{0}(\Omega)$ for every $k>1$, which forces $u$ to vanish on
the boundary of $\Omega$.
\end{uw}

In Section \ref{sec4} we give a different, analytic definition of
a solution of (\ref{eq1.1}) and we prove that actually it is
equivalent to the probabilistic definition. Before doing this we
would like to present the motivation behind the two definitions.
We begin with a concise presentation of famous Dynkin's formula.

Let $B\subset\Omega$ be a Borel set and let
\[
H^{1}_{0}(B)=\{u\in H^{1}_{0}(\Omega); u=0\mbox{ q.e. on }
B^{c}\}.
\]
(see \cite{Fuglede1,Fuglede2,KM}). Here and in the sequel for a
given function $u\in H^{1}_{0}(\Omega)$ we always consider its
quasi-continuous version. It is clear that $H^{1}_{0}(B)$ is a
closed subspace of the Hilbert space $H^{1}_{0}(\Omega)$.
Therefore
\[
H^{1}_{0}(\Omega)=H^{1}_{0}(B)\oplus \mathcal{H}_{\Omega\setminus
B},
\]
where $\mathcal{H}_{\Omega\setminus B}$ is the orthogonal
complement of $H^{1}_{0}(B)$. Let $H_{B}$ denote the operator of
the orthogonal projection on $\mathcal{H}_{\Omega\setminus B}$. By
\cite[Theorem 4.3.2]{Fukushima}, for every $u\in
H^{1}_{0}(\Omega)$,
\begin{equation}
\label{Projection}
H_{B}(u)(x)=E_{x}u(X_{\tau_{B}})\quad(=E_{x}u(X_{\tau_{B}\wedge\zeta}))
\end{equation}
for q.e. $x\in\Omega$, where
\[
\tau_{B}=\inf\{t>0,X_{t}\notin B\}.
\]

Let $G$ be a finely open subset of $\Omega$ and let
$(\EE_{G},D[\EE_{G}])$ denote the restriction of the form defined
by (\ref{DF}) to $G$, i.e.
\[
\EE_{G}(u,v)=\EE(u,v),\quad u,v\in D[\EE_{G}]=H^{1}_{0}(G).
\]
$(\EE_{G},D[\EE_{G}])$ is again a Dirichlet form (it may be no
longer regular). Let $\{p^{G}_{t},\,t\ge 0\}$ denote the
associated $C_{0}$-semigroup and $\{R^{G}_{\alpha},\, \alpha\ge
0\}$ the associated resolvent. By Theorems 4.4.2 and 4.4.4 in
\cite{Fukushima}, $(\EE_{G}, D[\EE_{G}])$ is transient and for
q.e. $x\in\Omega$,
\begin{equation}
\label{Parts} (R^{G}f)(x)=E_{x}\int_{0}^{\tau_{G}\wedge\zeta}
f(X_t)\,dt
\end{equation}
where $R^{G}=R^{G}_{0}$, which is well defined due to the
transiency of $(\EE_{G}, D[\EE_{G}])$. Dynkin's formula (see
\cite[page 153]{Fukushima}) says that for every finely open set
$G\subset\Omega$,
\[
(Rf)(x)=(R^{G}f)(x)+[H_{G}(Rf)](x)
\]
for q.e. $x\in\Omega$.

Suppose now that $\mu\in H^{-1}(\Omega)$ and there exists a weak
solution $u\in H^{1}_{0}(\Omega)$ of (\ref{eq1.1}) such that
$f(u)\in{L}^{2}(\Omega;m)$. Then for every $v\in
H^{1}_{0}(\Omega)$,
\begin{equation}
\label{eq2.13}
\EE(u^{k},v)=(f^{k}(u),v)_{L^2(\Omega;m)}+\langle\mu^{k},v\rangle,\quad
k=1,\dots,n
\end{equation}
that is
\[
u^{k}=Rf^{k}(u)+R\mu^{k},\quad k=1,\dots,n.
\]
Applying the operator $H_{G}$ to both sides of the above equation
and using Dynkin's formula we get
\begin{align*}
H_{G}(u^{k})=Rf^{k}(u)-R^{G}f^{k}(u)+R\mu^{k}-R^{G}\mu^{k}
=u^{k}-R^{G}f^{k}(u)-R^{G}\mu^{k}.
\end{align*}
As a consequence,
\begin{equation}
\label{Spectral} u^{k}=H_{G}(u^{k})+R^{G}
f^{k}(u)+R^{G}\mu^{k},\quad k=1,\dots,n.
\end{equation}
This equation expresses the property that if $u$ is a solution of
(\ref{eq1.1}) on $\Omega$ then for every finely open set $G$,
$u_{|G}$ is a solution of (\ref{eq1.1}) on $G$ with the boundary
condition $u_{|G}=u$ on $\partial G$.

In general, if $\mu\in\MM_{0,b}$, it is natural to look for
solutions of (\ref{eq1.1}) in the class of quasi-continuous
functions vanishing at the boundary of $\Omega$ and, roughly
speaking, such that they coincide with functions from the space
$H^1_0(\Omega)$ on each set $G$ from some family of finely open
set which covers $\Omega$. Therefore it is natural to require the
solution $u$ of (\ref{eq1.1}) to satisfy (\ref{Spectral}) or
(\ref{eq2.13}) for each set $G$ from some family of suitably
chosen (depending on $u$ in general) finely open sets. It is not
easy to deal with such families of equations. Fortunately, we can
obtain (\ref{Spectral}) from (\ref{eq.df.equiv}), and, in view of
Remark \ref{uw2.4}, from (\ref{eq.df}) if we know that $f(u)\in
qL^1(\Omega)$. Indeed, by standard arguments one can replace $T$
in (\ref{eq.df.equiv}) by $\tau_{G}$ with arbitrary finely open
set $G\subset\Omega$ and then putting $t=0$ and taking expectation
one can get
\[
u(x)=E_{x}u(X_{\tau_{G}\wedge\zeta})+E_{x}\int_{0}^{\tau_{G}\wedge\zeta}
f(u)(X_t)\,d\theta +E_{x}\int_{0}^{\tau_{G}\wedge\zeta}
dA^{\mu}_t,
\]
which in view of (\ref{Projection}) and (\ref{Parts}) gives
(\ref{Spectral}). The stochastic equations (\ref{eq.df}),
(\ref{eq.df.equiv}) are much more convenient to work with than
systems of the form (\ref{Spectral}). One of the major advantage
of (\ref{eq.df}) (resp. (\ref{eq.df.equiv})) lies in the fact that
it is well defined whenever $f\in q{L}^{1}_{loc}(\Omega)$ (resp.
$f(u)\in q{L}^{1}(\Omega)$) and $\mu\in \MM_{0,b}$. Moreover,
(\ref{eq.df}), (\ref{eq.df.equiv}) allow one to apply stochastic
analysis methods to study partial differential equations and are
well suited for dealing with the fine topology.

\nsubsection{Existence and uniqueness of probabilistic solutions}
\label{sec3}

We begin with the uniqueness result.
\begin{stw}
\label{uniqueness} Assume that \mbox{\rm(A4$'$)} is satisfied.
Then there exists  at most one probabilistic solution of
\mbox{\rm(\ref{eq1.1})}.
\end{stw}
\begin{dow}
Assume that $u_{1}, u_{2}$ are solutions of (\ref{eq1.1})
and $M_{1}, M_{2}$ are local MAFs associated with $u_{1},
u_{2}$, respectively. Then denoting $u=u_{1}-u_{2}$ and
$M=M_{1}-M_{2}$ we have
\[
u(X_{t})=u(X_{\tau})+\int_{t}^{\tau}(f(u_{1})-f(u_{2}))(X_{\theta})\,d\theta
-\int_{t}^{\tau}\,dM_{\theta}\quad 0\le t\le \tau,\quad
P_x\mbox{-a.s.}
\]
for q.e. $x\in\Omega$. By the It\^o-Tanaka formula and (A4$'$),
\begin{align*}
|u(X_{t})|&\le |u(X_{\tau})|+\int_{t}^{\tau} \langle
f(u_{1})-f(u_{2})(X_{\theta}),\hat{u}(X_{\theta})\rangle\,d\theta
-\int_{t}^{\tau}\hat{u}(X_{\theta})\,dM_{\theta}\\& \le
|u(X_{\tau})|-\int_{t}^{\tau}\hat{u}(X_{\theta}) \,dM_{\theta},
\quad 0\le t\le\tau, \quad P_x\mbox{-a.s.}
\end{align*}
for q.e. $x\in\Omega$. Without loss of generality we may assume
that $\int_{0}^{\cdot}\hat{u}(X_{\theta})\,dM_{\theta}$ is a true
martingale (otherwise one can apply the standard localization
procedure). Therefore putting $t=0$ and taking the expectation
with respect to $P_x$ we conclude that
\begin{equation}
\label{eq3.1} |u(x)|\le E_{x}|u(X_{\tau})|
\end{equation}
for q.e. $x\in\Omega$. Let $\{\tau_{k}\}$ be a sequence of
stopping times such that $0\le\tau_{k}<\zeta$, $k\ge 1$, and
$\tau_{k}\rightarrow \zeta$.  Since $u$ is of class (FD) and
$u\in\mathcal{C}_{0}(\Omega)$,  replacing $\tau$ by $\tau_{k}$ in
(\ref{eq3.1}) and then letting $k\rightarrow+\infty$ we conclude
that $|u|=0$ q.e.
\end{dow}

\begin{uw}
\label{uw.team} In general, the class $\mathcal{C}_{0}(\Omega)$ is
too large to ensure uniqueness of a solution of (\ref{eq1.1})
under (A4$'$). To see this, let us set $n=1$,
$\Omega=B(0,1)\equiv\{x\in\BRD;\, |x|<1\}$ and
$u(x)=\frac{1}{d-2}(\frac{1}{|x|^{d-2}}-1),\, d\ge 3$. Then $u\in
\mathcal{C}_{0}(\Omega)$ and from the Fukushima decomposition (see
\cite[Theorem 5.5.1]{Fukushima}) it follows  that
\[
u(X_{t})=u(X_{T\wedge\zeta})+\int_{t}^{T\wedge\zeta}dM_{\theta},\quad
t\ge0
\]
for some $M\in\MM_{loc}$. Thus, $u$ is a solution of (\ref{eq1.1})
with $f\equiv 0,\, \mu\equiv 0$. Obviously, the other solution of
the above equation is $v\equiv 0$. In fact, it is known (see
\cite{DMOP}) that $u$ is a renormalized solution of (\ref{eq1.1})
with $f\equiv0$ and $\mu=\sigma_{d-1}\delta_{0}$ (which is not a
smooth measure for $d\ge 2$), where $\sigma_{d-1}$ is the measure
of $\partial B(0,1)$.
\end{uw}

Let us recall that for a given additive functional $A$ of
$\mathbb{X}$ its energy is given by
\[
e(A)=\lim_{t\searrow 0}\frac{1}{2t}E_{m}A^{2}_{t}
\]
whenever the limit exists. It is known that for fixed regular
Dirichlet form $(\EE, D[\EE])$ and $u\in D[\EE]$ the additive
functional $A^{[u]}_{t}\equiv u(X_{t})-u(X_{0})$ admits the
so-called Fukushima decomposition, i.e. for q.e. $x\in\Omega$,
\[
A^{[u]}_{t}=M^{[u]}_{t}+N^{[u]}_{t},\quad t\ge 0,\quad
P_{x}\mbox{-a.s.},
\]
where $M^{[u]}$ is a martingale additive functional of $\BX$ of
finite energy and $N^{[u]}$ is a continuous additive functional of
$\BX$ of zero energy. This  decomposition is unique (see
\cite[Theorem 5.2.2]{Fukushima}).

\begin{lm}
\label{lm.mr} Assume that $u\in H^{1}_{0}(\Omega)$. Then for q.e.
$x\in\Omega$,
\begin{equation}
\label{eq3.2} P_{x}(\int_{0}^{T}|\nabla u(X_t)|^2\,dt<+\infty,\,
T\ge 0)=1
\end{equation}
and
\[
M_{t}^{[u]}=\int_{0}^{t}\nabla u(X_{\theta})\,dB_{\theta}, \quad
t\ge 0,\quad P_{x}\mbox{\rm-a.s.}
\]
\end{lm}
\begin{dow}
That $\nabla u$ satisfies (\ref{eq3.2}) follows immediately from
the fact that $\nabla u\in{L}^{2}(\Omega;m)$ (see (5.2.21) in
\cite{Fukushima}). Let $\{u_{n}\}\subset C_{0}^{\infty}(\Omega)$
be such that $u_{n}\rightarrow u$ strongly in $H_0^{1}(\Omega)$.
By \cite[Corollary 5.6.2]{Fukushima},
\[
M^{[u_{n}]}_{t}=\intot \nabla u_{n}(X_{\theta})\,dB_{\theta},
\quad t\ge 0,\quad P_{x}\mbox{-a.s.}
\]
for q.e. $x\in\Omega$. By \cite[Theorem 5.2.2]{Fukushima},
\[
e(M^{[u_{n}]}-M^{[u]})=\int_{\Omega}|\nabla(u_{n}-u)(x)|^{2}\,m(dx)
\rightarrow0.
\]
On the other hand, by (5.2.8) in \cite{Fukushima},
\[
e(M^{[u_{n}]}-\int_{0}^{\cdot}\nabla u(X_{\theta})\,dB_{\theta}
)=\int_{\Omega} |\nabla(u_{n}-u)(x)|^{2}\,m(dx)\rightarrow 0.
\]
Therefore the desired result follows from \cite[Theorem
5.2.1]{Fukushima}.
\end{dow}

\begin{lm}
\label{lm.aga} Let $B\subset\Omega$ be a Borel set such that
$m(\Omega\setminus B)=0$. Then for q.e. $x\in\Omega$,
\[
P_{x}(X_{t}\in B\mbox{ for a.e. }0\le t<\zeta )=1.
\]
\end{lm}
\begin{dow}
Let $A^{1}_{t}=\int_{0}^{t}\mathbf{1}_{B}(X_{\theta})\,d\theta,\,
A^{2}_{t}=\int_{0}^{t}\mathbf{1}_{\Omega}(X_{\theta})\,d\theta$.
Then  $A^{1}, A^{2}$ are PCAFs of $\BX$ and their associated Revuz
measures are $\mathbf{1}_{B}\cdot m$ and $m$, respectively. Since
$\mathbf{1}_{B}\cdot m=m$, it follows from uniqueness of the Revuz
correspondence that $A^{1}_{t}=A^{2}_{t}$, $t\ge 0$, $P_{x}$-a.s.
for q.e. $x\in\Omega$, which leads to the  desired result.
\end{dow}
\medskip

Let $\mathcal{FS}^{q}$, $q>0$, denote the set of all functions
$u\in \mathcal{C}(\Omega)$ such that for q.e. $x\in\Omega$,
\[
E_{x}\sup_{t\ge 0} |u(X_{t})|^{q}<+\infty,
\]
and let $\dot{H}^{1}_{loc}(\Omega)$ denote the space of all Borel
measurable functions on $\Omega$ for which there exists a
quasi-total family $\{U_{\alpha},\alpha\in I\}$ (i.e.
$\mbox{cap}(\Omega\setminus(\bigcup_{\alpha\in I} U_{\alpha}))=0$)
of finely open subsets of $\Omega$ such that for every $\alpha\in
I$ there exists a function $u_{\alpha}\in H^{1}_{0}(\Omega)$ such
that $u=u_{\alpha}$ q.e. on $U_{\alpha}$. For any function
$u\in\dot{H}^{1}_{loc}(\Omega)$ one can define its gradient as
\[
\nabla u=\nabla u_{\alpha}\mbox{ on }U_{\alpha}.
\]


\begin{tw}
\label{main} Assume that \mbox{\rm(A1)--(A4)} are satisfied. Then
there exists a solution $u$ of \mbox{\rm(\ref{eq1.1})} such that
$u\in \mathcal{FS}^{q}$ for $q\in (0,1)$, $u\in \dot{H}_{loc}^{1}(\Omega)$ and
\[
M_{t}=\int_{0}^{t} \nabla u(X_{\theta})\, dB_{\theta}, \quad t\ge
0,\quad P_{x}\mbox{\rm-a.s.}
\]
for q.e. $x\in \Omega$.
\end{tw}
\begin{dow}
Let $T_r$, $r>0$, denote the truncature operator, i.e.
\[
T_r(y)=\frac{ry}{|y|\vee r}\,,\quad y\in\BRD.
\]
Put $f_{n}=T_{n}(f)$, $n\in\BN$. Then $f_{n}$ is bounded and
satisfies (A4). Let $\{F_{n}\}$ be a generalized nest such that
$\mu_n=\mathbf{1}_{F_n}\cdot\mu\in H^{-1}(\Omega)$, $n\in\BN$ (for
the existence of such family see \cite[Theorem 2.2.4]{Fukushima}).
It is well known that there exists a solution $u_{n}\in
H^{1}_{0}(\Omega)$ of (\ref{eq1.1}) with $f_{n}$ in place of $f$
and $\mu_{n}$ in place of $\mu$. By \cite{KR}, $u_{n}\in
\mathcal{FS}^{2}$ and there exists $M^{n}\in\MM^{2}$ such that
for q.e. $x\in\Omega$,
\begin{align*}
u_{n}(X_{t})&=u_{n}(X_{T\wedge\zeta})+\int_{t}^{T\wedge\zeta}
f^{k}_{n}(u_{n})(X_{\theta})\,d\theta+\int_{t}^{T\wedge\zeta}
dA^{\mu_{n}}_{\theta} \\
&\quad-\int_{t}^{T\wedge\zeta} dM^{n}_{\theta},\quad 0\le t\le
T<\infty,\quad P_{x}\mbox{-a.s.}
\end{align*}
As in \cite[page 201]{Fukushima}) one can check that the CAFs
$\int_{0}^{\cdot}f^{k}_{n}(u_{n})(X_{\theta})\,d\theta$,
$A^{\mu_{n}}_{t}$ are of zero energy. It follows from uniqueness
of the Fukushima decomposition and Lemma \ref{lm.mr} that for q.e.
$x\in \Omega$,
\[
M^{n}_{t}=\int_{0}^{t}\nabla u_{n}(X_{\theta})\,dB_{\theta}, \quad
t\ge 0,\quad P_{x}\mbox{-a.s.}
\]
Since $u_{n}$ is of class (FD) and $u\in\mathcal{C}_{0}(\Omega)$,
in much the same way as in the proof of (\ref{eq2.8}) we show that
\begin{equation}
\label{eq2.1} |u_{n}(x)|\le E_{x}\int_{0}^{\zeta}d|A^{\mu}|_t
\end{equation}
for q.e. $x\in\Omega$. Write
\[
v(x)= E_{x}\int_{0}^{\zeta}d|A^{\mu}|_t.
\]
By \cite{KR}, $v\in \FF\mathcal{S}^{q}$, $q\in (0,1)$, for every
$k>0$, $T_{k}(v)\in H^{1}_{0}(\Omega)$ and $v$ is of class (FD).
Let $G_{k}=\{|v|<k\}$. Since $v$ is quasi-continuous, $G_{k}$ is
finely open. Moreover, the family $\{G_{k}\}$ is quasi-total. Let
us put
\[
\tau_{k}=\inf\{t>0; X_{t}\in G_{k}\}.
\]
By It\^o's formula, for q.e. $x\in\Omega$ we have
\begin{align*}
&|u_{n}(x)|^2+E_{x}\int_{0}^{\zeta\wedge\tau_{k}}|\nabla
u_{n}|^{2}(X_t)\,dt\\
&\qquad =E_{x}|u_{n}(X_{\tau_{k}\wedge\zeta})|^2
+2E_{x}\int_{0}^{\zeta\wedge \tau_{k}}\langle
f^{n}(u_{n})(X_t), u_{n}(X_t)\rangle\,dt\\
&\qquad\quad +2E_{x}\int_{0}^{\zeta\wedge\tau_{k}}\langle
u_{n}(X_t),dA^{\mu}_t\rangle.
\end{align*}
By the definition of $\tau_{k}$, (\ref{eq2.1}) and (A4), for q.e.
$x\in\Omega$,
\begin{equation}
\label{eq3.06} E_{x}\int_{0}^{\zeta\wedge\tau_{k}}|\nabla
u_{n}(X_t)|^2\,dt \le
2k+2kE_{x}\int_{0}^{\zeta\wedge\tau_{k}}d|A^{\mu}|_t.
\end{equation}
Since $\Omega$ is bounded, there exists $R>0$ such that
$\Omega\subset B(0,R)$. Hence
\begin{equation}
\label{eq.esp} E_{x}\zeta\le E_{x}\tau_{B(0,R)}\le
C(d)(R^{2}-|x|^{2})
\end{equation}
for every $x\in\Omega$ (for the last inequality see, e.g.,
\cite[page 253]{KaratzasShreve}). From (\ref{RD1}) and
(\ref{eq.esp}) it follows in particular that for every
$f\in\mathcal{B}^{+}(\Omega)$ and $\mu\in S$,
\begin{equation}
\label{TV} E_{m}\int_{0}^{\zeta}f(X_t)\,dA^{\mu}_t\le C(\Omega,d)
\|f\cdot\mu\|_{TV},
\end{equation}
where $\|\cdot\|_{TV}$ denotes the total variation norm. By
(\ref{eq3.06}) and (\ref{TV}),
\begin{equation}
\label{eq2.2} E_{m}\int_{0}^{\zeta\wedge\tau_{k}}|\nabla
u_{n}(X_t)|^2\,dt\le 2km(\Omega)+2kc|\mu|(\Omega).
\end{equation}
Since $G_{k}$ is finely open, it follows from  \cite[Theorem
4.2.2]{Fukushima} that $\EE_{G_{k}}$ is a regular Dirichlet form
on ${L}^{2}(G_{k};m)$ and the semigroup $\{p^{G_{k}}_{t},t\ge 0\}$
is determined by the process $\mathbb{X}^{G_{k}}$ (see
(\ref{Parts})). Therefore for every $f\in {L}^{2}(G_{k};m)$,
\[
R^{G_{k}}f(x)=E_{x}\int_{0}^{\tau_{k}\wedge\zeta}f(X_{t})\,dt
\mbox{ for }m\mbox{-a.e. }x\in G_{k}.
\]
Moreover, by \cite[Theorem 4.4.4]{Fukushima}, $\EE_{G_{k}}$ is
transient and $D[\EE_{G_{k}}]=H^{1}_{0}(G_{k})$. Therefore from
(\ref{eq2.2}) it follows that
\begin{equation}
\label{eq.res} \sup_{n\ge 1}\|R^{G_{k}}(|\nabla
u_{n}|^{2})\|_{{L}^{1}(\Omega;m)}<+\infty.
\end{equation}
On the other hand, by \cite[Lemma 5.1.10]{Fukushima},
\[
\|R^{G_{k}}(|\nabla u_{n}|^{2})\|_{{L^1}(\Omega;m)}
=\int_{G_{k}}|\nabla u_{n}|^{2}(y)R^{G_{k}}1(y)\,m(dy).
\]
Since $R^{G_{k}}1\in D[\EE_{G_{k}}]=H^{1}_{0}(G_{k})\subset
H^{1}_{0}(\Omega)$ and $|u_{n}(x)|\le k$ on $G_{k}$ for q.e.
$x\in\Omega$, we conclude from the above estimate that
\begin{equation}
\label{UE} \sup_{n\ge 1}\int_{\Omega}|\nabla(u_{n}\cdot
R^{G_{k}}1)|^{2}(y)\,dy<+\infty.
\end{equation}
Since ${L}^{2}(\Omega;m)$ has the Banach-Saks property, it follows
from (\ref{UE}) that one can choose a subsequence (still denoted
by $\{n\}$) such that $\sigma_{n}(\{\nabla(u_{n}\cdot
R^{G_{k}}1)\})$ is convergent in ${L}^{2}(\Omega;m)$ for every
$k\ge 1$. By \cite[Theorem 2.1.4]{Fukushima}, one can find a
further subsequence (still denoted by $\{n\}$) such that
$\{\sigma_{n}(\{u_{n}R^{G_{k}}1\})\}$ is convergent q.e. for every
$k\ge 1$. Since $G_{k}$ is finely open,
$R^{G_{k}}1(x)=E_{x}(\tau_{k}\wedge\zeta)>0$ q.e. on $G_{k}$.
Therefore $\{\sigma_{n}(\{u_{n}\})\}$ is convergent q.e. on
$\Omega$. Set $u(s,x)=\limsup_{n\rightarrow
+\infty}\sigma_{n}(\{u_{n}\})(s,x)$ for $(s,x)\in\Omega$. Since
$R^{G_{k}}1\in H^{1}_{0}(\Omega)$, $u$  is quasi-continuous. Using
this and the fact that $R^{G_{k}}1>0$ q.e. on $G_{k}$ we see that
one can find a quasi-total family $\{\tilde{G}_{k}\}$ such that
for every $k\ge1$,
\[
\int_{\tilde{G}_{k}}|\sigma_{n}(\{\nabla u_{n}\})
-\sigma_{m}(\{\nabla u_{m}\})|^{2}(y)\,dy\rightarrow 0
\]
as $n\rightarrow+\infty$. Therefore we may define a measurable
function $w:\Omega\rightarrow\mathbb{R}^{n}\times\BRD$ such that
$w_{|\tilde{G}_{k}}=\lim_{n\rightarrow +\infty}\sigma_{n}(\{\nabla
u_{n}\})$ in ${L}^{2}(\tilde{G}_{k};m)$ for every $k\ge 1$. Let us
fix $\alpha>0$ and $\nu\in S^{(0)}_{00}$ (see \cite[Section 2]{Fukushima}). Then
\begin{align*}
P_{\nu}(\int_{0}^{\tilde{\tau}_{k}}|w(X_t)|^2\,dt>\alpha) &\le
\alpha^{-1}E_{\nu}
\int_{0}^{\tilde{\tau}_{k}}\mathbf{1}_{\tilde{G}_{k}}|w(X_t)|^2\,
dt\\
&=\alpha^{-1}\langle \nu,
R^{\tilde{G}_{k}}(\mathbf{1}_{\tilde{G}_{k}}w^{2})\rangle\\
&=\alpha^{-1}\langle R^{\tilde{G}_{k}}\nu,
\mathbf{1}_{\tilde{G}_{k}}w^{2}\rangle\\
&\le\alpha^{-1}\|R\nu\|_{\infty}\int_{\tilde{G}_{k}}|w(y)|^{2}\,dy.
\end{align*}
Similarly,
\begin{align*}
P_{\nu}(\int_{0}^{\tilde{\tau_{k}}}|(\sigma_{n}(\{\nabla
u_{n}\})-w)(X_t)|^{2}\,dt>\alpha) \le
\alpha^{-1}\|R\nu\|_{\infty}\int_{\tilde{G}_{k}}|\sigma_{n}(\{\nabla
u_{n}\})(y)-w(y)|^{2}\,dy,
\end{align*}
which converges to 0 as $n\rightarrow\infty$. Using the above two
inequalities and the Borel-Cantelli lemma one can show that
$P_{x}(\int_{0}^{T}|w(X_t)|^2\,dt<+\infty,\, T\ge 0)=1$ for q.e.
$x\in\Omega$ and there exists a subsequence (still denoted by
$\{n\}$) such that for every $T>0$,
\begin{equation}
\label{eq3.4} \int_{0}^{\cdot} \sigma_{n}(\{\nabla
u_{n}\})(X_t)\,dB_t \rightarrow \int_{0}^{\cdot} w(X_t)\,dB_t
\end{equation}
in ucp on $[0,T]$ with respect to $P_{x}$ for q.e. $x\in \Omega$
(see the proof of \cite[Theorem 5.2.1]{Fukushima}). Furthermore,
by (\ref{UE}), the Rellich-Kondrachov theorem and the fact that
$R^{G_{k}}1>0$ q.e. on $G_{k}$ and $\{G_{k}\}$ is a quasi-total
family we conclude that there exists a subsequence (still denoted
by $\{n\}$) such that $u_{n}\rightarrow u$ $m$-a.e. Hence, by
Lemma \ref{lm.aga}, for $n\ge k$ and $0\le\tau<\zeta$ we have
\begin{align*}
&P_{x}(\int_{0}^{\tau_{k}\wedge\tau}|(f^{n}(u_{n})-f(u))(X_t)|\,
dt>\varepsilon)\\
&\qquad=P_{x}(\int_{0}^{\tau_{k}\wedge\tau}|(f(u_{n})-f(u))(X_t)|\,
dt>\varepsilon)\\
&\qquad=P_{x}(\int_{0}^{\tau_{k}\wedge\tau}\mathbf{1}_{\{|u_{n}|\le
k\}}\mathbf{1}_{\{|u|\le k\}}
|(f(u_{n})-f(u))(X_t)|\,dt>\varepsilon),
\end{align*}
which for q.e. $x\in\Omega$ converges to 0 as
$n\rightarrow+\infty$. From this we conclude that for q.e. $x\in
\Omega$, under $P_x$,
\begin{equation}
\label{eq3.5}
\int_{0}^{\cdot\wedge\tau}|(f_{n}(u_{n})-f(u))(X_t)|\,dt\rightarrow0
\end{equation}
in ucp on $[0,T]$ for every $0\le\tau<\zeta$. Moreover,
\begin{align}
\label{eq3.6}
\sigma_{n}(\{u_{n}\})(X_{t})&=\sigma_{n}(\{u_{n}\})(X_{T\wedge\tau})
+\int_{t}^{T\wedge\tau}\sigma_{n}(\{f^{n}(u_{n})\})(X_{\theta})\,d\theta
\nonumber\\
&\quad+\int_{t}^{T\wedge\tau}\sigma_{n}(\{dA^{\mu_{n}}_{\theta}\})
-\int_{t}^{T\wedge\tau}\sigma_{n}(\{\nabla
u_{n}\})(X_{\theta})\,dB_{\theta}
\end{align}
and without loss of generality (see \cite[Theorem
4.1.1]{Fukushima}) we may assume that $\sigma_{n}(\{u_{n}\})$ is
convergent on $\Omega$ except for some properly exceptional set
$N\subset \Omega$ (see (\ref{exceptional})). Therefore letting
$n\rightarrow\infty$ in (\ref{eq3.6}) and using (\ref{eq3.4}),
(\ref{eq3.5}) we see that for every $T>0$,
\begin{align*}
u(X_{t})&=u(X_{T\wedge\tau})+\int_{t}^{T\wedge\tau}f(u)(X_{\theta})\,d\theta
+\int_{t}^{T\wedge\tau}dA^{\mu}_{\theta}\\
&\quad -\int_{t}^{T\wedge\tau}w(X_{\theta})\,dB_{\theta},\quad
0\le t\le T\wedge\tau,\quad P_x\mbox{-a.s.}
\end{align*}
for q.e. $x\in \Omega$. From this we get (\ref{eq.df}). From
(\ref{eq.df}), the fact that $|u(X_{t})|\le |v(X_{t})|$, $t\ge 0$,
$P_{x}$-a.s. for q.e. $x\in\Omega$ and Remark \ref{van.b} it
follows  that $u$ is a solution of (\ref{eq1.1}), $u\in
\mathcal{FS}^{q}$ for $q\in (0,1)$ and $u$ is of class (FD). By
(\ref{eq.res}),  $\{|\nabla u_{n}|^{2} T_{r}(R^{G_{k}}1)\}_n$ is
bounded in ${L}^{1}(\Omega;m)$ for every $r>0$. Since
$R^{G_{k}}1\in D[\EE_{G_{k}}]$, $T_{r}(R^{G_{k}}1)\in
D[\EE_{G_{k}}]$. From this we conclude that
\[
\sup_{n\ge 1}\int_{\Omega}
|\nabla(u_{n}T_{r}(R^{G_{k}}1))|^{2}(y)\,dy<+\infty,
\]
which implies that $uT_{r}(R^{G_{k}}1)\in H^{1}_{0}(\Omega)$ for
every $r>0$ and $k\ge 0$. Since $R^{G_{k}}1$ is quasi continuous
and positive q.e. on $G_{k}$, $\{R^{G_{k}}1>r, k\ge 1, r>0\}$
forms a finely open quasi-total family. This shows that
$u\in\dot{H}^{1}_{loc}(\Omega)$ and $w=\nabla u$ since
$r^{-1}uT_{r}(R^{G_{k}}1)=u$ q.e. on $\{R^{G_{k}}>r\}$.
\end{dow}

\nsubsection{Analytic solutions}
\label{sec4}

To formulate the definition of a solution of (\ref{eq1.1}) in the
analytic sense we will need the following lemma.

\begin{lm}
\label{lm2} Let $\mu$ be a smooth measure on $\Omega$. Then there
exists a finely open quasi-total family $\{G_{k}\}$ such that
$\mathbf{1}_{G_{k}}\mu\in H^{-1}(\Omega)$ for every $k\ge 0$.
\end{lm}

\begin{dow}
The proof is a slight modification of the proof of \cite[Lemma
5.1.7]{Fukushima}. Let $f$ be a Borel bounded positive function
and let
\[
\phi(x)=E_x\int_{0}^{\infty}e^{-t-A_{t}} f(X_{t})\,dt,
\]
where $A=A^{\mu}$. Then $\phi=R^{A}_{1}f$, where
$\{R^{A}_{\alpha}\}$ is the resolvent associated with the
perturbed form (see \cite{Fukushima}). It follows in  particular
that $\phi$ is quasi-continuous. Put
\[
G_{k}=\{x\in\Omega; \phi(x)>k^{-1}\},\quad
\mu_{k}=\mathbf{1}_{G_{k}}\cdot\mu.
\]
Since $f$ is positive and $\phi$ is quasi-continuous, $\{G_{k}\}$
is a finely open quasi-total family. We have
\[
R_{1}\mu_{k}(x)=E_{x}\int_{0}^{\zeta}
e^{-t}\mathbf{1}_{G_{k}}(X_{t}) \,dA_{t} \le
kE_{x}\int_{0}^{\zeta} e^{-t}\phi(x)\,dA_{t}\le kR^{1}f(x),
\]
the last inequality being a consequence of \cite[Lemma
5.1.5]{Fukushima}. Since $R^{1}f\in H^{1}_{0}(\Omega)$, it follows
from Theorem 2.2.1 and  Lemma 2.3.2 in \cite{Fukushima} that
$R_{1}\mu_{k}\in H^{1}_{0}(\Omega)$, hence that $\mu_{k}\in
H^{-1}(\Omega)$.
\end{dow}

\begin{uw}
\label{uwrepeat} Suppose that $\{U_{\alpha},\,\alpha\in T_{1}\}$,
$\{V_{\beta},\,\beta\in T_{2}\}$ are two finely open quasi-total
families. Then there exists a finely open quasi-total family
$\{W_{\gamma},\,\gamma\in T_{3}\}$ such that for every $\gamma\in
T_{3}$ there exists $\alpha\in T_{1}, \beta\in T_{2}$ such that
$W_{\gamma}\subset U_{\alpha}\cap V_{\beta}$. To see this, it
suffices to take $T_{3}=T_{1}\times T_{2}$ and
$\{W_{(\alpha,\beta)}=U_{\alpha}\cap V_{\beta},\,\alpha\in
T_{1},\,\beta_{1}\in T_{2}\}$.
\end{uw}

From the the definition of the space $\dot{H}^{1}_{loc}(\Omega)$,
Lemma \ref{lm2} and Remark \ref{uwrepeat} it follows that for
given $u\in\dot{H}^{1}_{loc}(\Omega)$, $f\in
q{L}^{1}_{loc}(\Omega)$ and $\mu\in S$ there exists a finely open
quasi-total family $\{U_{\alpha}, \alpha\in T\}$ such that for
every $\alpha\in T$,
\begin{equation}
\label{qT} u=u_{\alpha}\mbox{ q.e. on }U_{\alpha}\mbox{ for some }
u_{\alpha}\in H^{1}_{0}(\Omega), \quad
\mathbf{1}_{U_{\alpha}}f\in{L}^{2}(\Omega;m),\quad
\mathbf{1}_{U_{\alpha}}\cdot\mu\in H^{-1}.
\end{equation}
In what follows we say that $\{U_{\alpha},\,\alpha\in T\}$ is a
finely open quasi-total family for the triple
$(u,f,\mu)\in\dot{H}^{1}_{loc}(\Omega)\times
q{L}^{1}_{loc}(\Omega)\times S$ if (\ref{qT}) is satisfied for
every $\alpha\in T$.

Proposition \ref{uniqueness} and Remark \ref{uw.team} suggest that
the class $\dot{H}^{1}_{loc}(\Omega)\cap\mathcal{C}_{0}(\Omega)$
is too large to get uniqueness of analytic solutions of
(\ref{eq1.1}) under (A4$'$), and secondly, that the uniqueness
holds if we restrict the class of solutions to functions which
additionally are of class (FD). Unfortunately, we do not know how
to define the class (FD) analytically. Therefore to state the
definition of a solution of (\ref{eq1.1}) in a purely analytic way
we introduce a class $\mathcal{U}$, which is a bit narrower than
(FD).

$\mathcal{U}=\{u\in\mathcal{C}_{0}(\Omega): |u|\le v$ q.e., where
$v$ is a  solution in the sense of duality (see
\cite{Stampacchia}) of the problem
\begin{equation}
\label{eq4.7}
-\frac12\Delta v=\nu\mbox{ in }\Omega, \quad
v=0\mbox{ on }\partial\Omega
\end{equation}
for some $\nu\in\MM_{0,b}$\}.

\begin{df}
We say that $u:\Omega\rightarrow\mathbb{R}^{n}$ is a solution of
(\ref{eq1.1}) in the analytic sense if
\begin{enumerate}
\item[({\rm a})]$u\in \mathcal{U}$,
\item[({\rm b})]$u\in\dot{H}^{1}_{loc}(\Omega)\cap \mathcal{C}_{0}(\Omega)$,
$f(u)\in qL^{1}_{loc}(\Omega)$,
\item [({\rm c})]For some finely open quasi-total
family  $\{G_{l},l\ge 0\}$ for $(u,f(u),\mu)$,
\begin{equation}
\label{eq.df.a} \EE(u^{k},v)=(f^{k}(u),v)_{{L}^{2}(\Omega;m)}
+\langle \mu,v\rangle, \quad v\in H^{1}_{0}(G_{l})
\end{equation}
for every $k=1,\dots,n$ and $l\ge0$.
\end{enumerate}
\end{df}

\begin{prz}
Let $\Omega$, $u$ be as in Remark \ref{uw.team} and let
$G_{l}=\{x\in B(0,1); l^{-1}<|x|<1\}$. Then $\{G_{l},\,l\ge1\}$ is
a finely open quasi-total family for $u$,
$u\in\dot{H}^{1}_{loc}(\Omega)\cap\mathcal{C}_{0}(\Omega)$ and
$\EE(u^{k},v)=0$ for $v\in H^{1}_{0}(G_{l})$, i.e. $u$ satisfies
(b), (c) of the above definition with $f\equiv 0$, $\mu\equiv 0$.
Of course, $v\equiv 0$ satisfies the same conditions, too.
\end{prz}

\begin{lm}
\label{lm.sup} If $u\in H^{1}_{0}(\Omega)$ then $u$ is of class
\mbox{\rm(FD)}.
\end{lm}
\begin{dow}
If $u\in H^{1}_{0}(\Omega)$ then $u^{+}\in H^{1}_{0}(\Omega)$, so
we may assume that $u\ge 0$. Let $v\in H^{1}_{0}(\Omega)$ be such
that $v\ge u$ q.e. and $-\Delta v\ge 0$ in the distributional
sense (one can, for example, take $v$ to be a solution of the
obstacle problem with barrier $u$). Then there exists a
nonnegative diffuse measure $\mu\in H^{-1}(\Omega)$ such that
$-\frac12\Delta v=\mu$ in $H^{-1}(\Omega)$. By \cite{KR}, $v$ is
of class (FD), and hence so is $u$.
\end{dow}

\begin{stw}
\label{stw.12} Assume that $u\in\dot{H}^{1}_{loc}(\Omega)$ and
$f(u)\in qL^{1}_{loc}$. Then $u$  satisfies
\mbox{\rm{(\ref{eq.df})}} iff $u$ satisfies
\mbox{\rm{(\ref{eq.df.a})}}.
\end{stw}
\begin{dow}
Assume that $u$ satisfies (\ref{eq.df}). Let $\{G_{l},\,l\ge1\}$
be a finely open quasi-total family for $(u,f(u),\mu)$. Fix
$l\in\mathbb{N}$, $v\in H^{1}_{0}(G_{l})$ and let
$\tau_{l}=\tau_{G_{l}}$. By (\ref{eq.df}), for q.e. $x\in\Omega$
and every $0\le\tau<\zeta$ we have
\begin{align}
\label{star1}
u(X_{t})&=u(X_{\tau_{l}\wedge\tau})
+\int_{t}^{\tau_{l}\wedge\tau}f(u)(X_{\theta})\,d\theta
+\int_{t}^{\tau_{l}\wedge\tau} dA^{\mu}_{\theta} \nonumber\\
&\quad-\int_{t}^{\tau_{l}\wedge\tau}
dM_{\theta},\quad 0\le t\le \tau_{l}\wedge\tau,\quad P_x\mbox{-a.s.}
\end{align}
Let $u_{l}\in H^{1}_{0}(\Omega)$ be such that such that $u=u_{l}$
q.e. on $G_{l}$. Then by (\ref{star1}),
\begin{equation}
\label{eq4.2} u_{l}(x)=u_{l}(X_{\tau_{l}\wedge\tau})
+\int_{0}^{\tau_{l}\wedge\tau}f(u)(X_t)\,dt
+\int_{0}^{\tau_{l}\wedge\tau} dA^{\mu}_t
-\int_{0}^{\tau_{l}\wedge\tau} dM_t.
\end{equation}
By the definition of $\tau_{l}$ and the remark following
(\ref{RD2}),
\begin{align*}
E_{x}\int_{0}^{\tau_{l}\wedge\tau} |f(u)(X_t)|\,dt
&=E_{x}\int_{0}^{\tau_{l}\wedge\tau}\mathbf{1}_{G_{l}}(X_t)
|f(u)(X_t)|\,dt \\
&\le E_{x}\int_{0}^{\zeta}
\mathbf{1}_{G_{l}}|f(u)(X_t)|\,dt<+\infty,
\end{align*}
for q.e. $x\in\Omega$. By Fatou's lemma,
$E_{x}\int_{0}^{\tau_{l}\wedge\zeta} |f(u)(X_t)|\,dt<+\infty$ for
q.e. $x\in\Omega$. Likewise, by the remark following (\ref{RD2}),
$E_{x}\int_{0}^{\zeta} d|A^{\mu}|_t<+\infty$ and, by \cite[Theorem
4.3.2]{Fukushima}, $E_{x}|u(X_{\tau_{l}\wedge\zeta})|<+\infty$ for
q.e. $x\in\Omega$. Therefore from (\ref{eq4.2}), the fact that
$u_{l}\in\mathcal{C}_{0}(\Omega)$ (see Remark \ref{van.b}) and
$u_{l}$ is of class (FD) (see Lemma \ref{lm.sup}) it follows that
for q.e. $x\in\Omega$,
\begin{equation}
\label{star2} u_{l}(x)=E_{x}u_{l}(X_{\tau_{l}\wedge\zeta})
+E_{x}\int_{0}^{\tau_{l}\wedge\zeta}f(u)(X_t)\,dt
+E_{x}\int_{0}^{\tau_{l}\wedge\zeta}dA^{\mu}_t.
\end{equation}
By (\ref{Projection}) and (\ref{Parts}) one can rewrite the above
equation in the form
\begin{equation}
\label{a}
u_{l}(x)=(H_{G_{l}}u_{l})(x)+(R^{G_{l}}f(u))(x)+(R^{G_{l}}\mu_{l})(x).
\end{equation}
Hence
\begin{equation}
\label{b}
(u_{l},v)_{{L}^{2}(\Omega;m)}=(H_{G_{l}}(u_{l}),v)_{{L}^{2}(\Omega;m)}
+(R^{G_{l}}f,v)_{{L}^{2}(\Omega;m)}
+(R^{G_{l}}\mu_{l},v)_{{L}^{2}(\Omega;m)},
\end{equation}
that is
\begin{equation}
\label{c}
(u_{l}-H_{G_{l}}(u_{l}),v)_{{L}^{2}(\Omega;m)}=(R^{G_{l}}f,v)_{{L}^{2}(\Omega;m)}
+(R^{G_{l}}\mu_{l},v)_{{L}^{2}(\Omega;m)}.
\end{equation}
Since $v, u^{k}_{l}-H_{G_{l}}(u^{k}_{l})\in
H^{1}_{0}(G_{l})\subset H^{1}_{0}(\Omega)$,
\[
(u^{k}_{l}-H_{G_{l}}(u^{k}_{l}), v)_{{L}^{2}(\Omega;m)}
=\EE(u^{k}_{l}-H_{G_{l}}(u^{k}_{l}), R^{G_{l}}v),\quad
k=1,\dots,n,
\]
and hence
\[
\EE(u^{k}_{l}-H_{G_{l}}(u^{k}_{l}), R^{G_{l}}v)=\EE(u^{k}_{l},
R^{G_{l}}v),\quad k=1,\dots,n,
\]
because $H^{1}_{0}(G_{l})$ and $\HH_{\Omega\setminus G_{l}}$ are
orthogonal and  $R^{G_{l}}v\in H^{1}_{0}(G_{l})$. Using this and
symmetry of $R^{G_{l}}$ we conclude from (\ref{c}) that
\[
\EE(u^{k}_{l}, R^{G_{l}}v)=(f^{k}(u),
R^{G_{l}}v)_{{L}^{2}(\Omega;m)}
+(\mu_{l}^{k},R^{G_{l}}v)_{{L}^{2}(\Omega;m)},\quad k=1,\dots,n,
\]
hence that
\begin{equation}
\label{d} \EE(u^{k}, R^{G_{l}}v)=(f^{k}(u),
R^{G_{l}}v)_{{L}^{2}(\Omega;m)}
+(\mu^{k},R^{G_{l}}v)_{{L}^{2}(\Omega;m)},\quad k=1,\dots,n
\end{equation}
since $R^{G_{l}}v\in H^{1}_{0}(G_{l})$. Finally, since (\ref{d})
holds true for every $v\in H^{1}_{0}(G_{l})$, $l\in \mathbb{N}$,
and $R^{G_{l}}(H^{1}_{0}(G_{l}))$ is a dense subset of
$H^{1}_{0}(G_{l})$ in the  norm induced by the form $\EE_{G_{l}}$,
i.e. the norm defined as
$\|v\|_{\EE_{G_{l}}}=\EE_{G_{l}}(v,v)+(v,v)_{{L}^{2}(\Omega,m)}$
for $v\in H^{1}_{0}(G_{l})$, we conclude from (\ref{d})  that for
every $v\in H^{1}_{0}(G_{l})$,
\[
\EE(u^{k}, v)=(f^{k}(u),v)_{{L}^{2}(\Omega;m)}
+(\mu^{k},v)_{{L}^{2}(\Omega;m)},\quad k=1,\dots,n,
\]
which shows that $u$ satisfies (\ref{eq.df}).

Now, let us assume that $u$ satisfies (\ref{eq.df.a}).
Let $\{G_{l},\, l\ge 1\}$ be a finely open quasi-total family
for $(u,f(u),\mu)$. Then taking $v=R^{G_{l}}\nu$ with
$\nu\in\MM_{0,b}\cap H^{-1}(\Omega)$ as a test function in
(\ref{eq.df.a}) and reversing the steps of the proof of (\ref{d})
we get (\ref{b}) with $v$ replaced by $\nu$ and with the duality
$\langle\cdot,\cdot\rangle$ in place of
$(\cdot,\cdot)_{{L}^{2}(\Omega;m)}$, i.e.
\[
\langle u_{l},\nu\rangle=\langle H_{G_{l}}(u_{l}),\nu\rangle
+\langle R^{G_{l}}f,\nu\rangle +\langle
R^{G_{l}}\mu_{l},\nu\rangle.
\]
By \cite[Theorem 2.2.3]{Fukushima}, this implies (\ref{a}), which
in turn gives (\ref{star2}). Now we show that (\ref{star2})
implies (\ref{star1}). To this end, let us first note that by
\cite[Theorem 4.1.1]{Fukushima} one can assume that (\ref{star2})
holds except for some properly exceptional set. Let
$\tau\in\mathcal{T}$ be such that $0\le\tau <\tau_{l}\wedge\zeta$.
Then
\begin{align*}
u_{l}(X_{\tau})=E_{X_{\tau}}u_{l}(X_{\tau_{l}\wedge\zeta})
+E_{X_{\tau}}\int_{0}^{\tau_{l}\wedge\zeta}f(u)(X_{t})\,dt
+E_{X_{\tau}}\int_{0}^{\tau_{l}\wedge\zeta}dA^{\mu}_{t},\quad
P_{x}\mbox{-a.s.}
\end{align*}
Let $\{\theta_{t}, \,t\ge0\}$ denote the family of shift
operators. By the strong Markov property and the fact that
$(\tau_{l}\wedge\zeta)\circ\theta_{\tau}-\tau=\tau_{l}\wedge\zeta$,
\begin{align*}
u_{l}(X_{\tau})=E_{x}\left( u_{l}(X_{\tau_{l}\wedge\zeta})
+\int_{\tau}^{\tau_{l}\wedge\zeta}f(u)(X_{t})\,dt
+\int_{\tau}^{\tau_{l}\wedge\zeta}dA^{\mu}_{t}|\FF_{\tau}\right),
\quad P_{x}\mbox{-a.s.}
\end{align*}
Write
\[
M_{t}^{l,x}=E_{x}\left(
\int_{0}^{\tau_{l}\wedge\zeta}f(u)(X_{\theta})\,dt
+\int_{0}^{\tau_{l}\wedge\zeta}
dA^{\mu}_{\theta}|\FF_{t}\right)-u_{l}(X_{0}), \quad t\ge 0.
\]
By \cite[Lemma A.3.5]{Fukushima}, $M^{l,x}$ has a continuous
version $M^{l}$ which does not depend on $x$. Observe that $M^{l}$
is a MAF of $\BX$ and
\[
u_{l}(X_{\tau})=u(X_{0})-\int_{0}^{\tau}f(u)(X_{t})\,dt
-\int_{0}^{\tau}dA^{\mu}_{t}+\int_{0}^{\tau} dM^{l}_{t}, \quad
P_{x}\mbox{-a.s.}
\]
Since the above equality is satisfied for every $\tau\in
\mathcal{T}$ such that $0\le\tau<\tau_{l}\wedge\zeta$, we can
assert that
\[
u_{l}(X_{t})=u(X_{0})-\int_{0}^{t}f(u)(X_{\theta})\,d\theta
-\int_{0}^{t} dA^{\mu}_{\theta}+\int_{0}^{t}dM^{l}_{\theta}, \quad
0\le t\le\tau_{l}\wedge\zeta,\quad P_{x}\mbox{-a.s.}
\]
Without loss of generality we may assume that the family
$\{G_{l},\,l\ge 1\}$ is ascending. Then by \cite[Lemma
5.5.1]{Fukushima}, $M^{l}_{t}=M^{l+1}_{t}$, $0\le t\le \tau_{l}$.
Therefore the process $M$ defined as $M_{t}=\lim_{l\rightarrow
+\infty} M^{l}_{t}$, $t\ge 0$, is well defined, $M$ is a local MAF
of $\BX$ and
\[
u_{l}(X_{t})=u(X_{0})-\int_{0}^{t}f(u)(X_{\theta})\,d\theta
-\int_{0}^{t} dA^{\mu}_{\theta}+\int_{0}^{t} dM_{\theta}, \quad
0\le t\le\tau_{l}\wedge\zeta,\quad P_{x}\mbox{-a.s.},
\]
which implies (\ref{star1}). Finally, letting $l\rightarrow
+\infty$ in (\ref{star1}) we get (\ref{eq.df}).
\end{dow}

\begin{stw}
\label{stw4.6} If $u\in \mathcal{U}$ then $u$ is of class
\mbox{\rm(FD)}.
\end{stw}
\begin{dow}
From \cite{KR} it follows that if $\nu\in\MM_{0,b}$ then the
solution of (\ref{eq4.7}) in the sense of duality  is of class
(FD). This implies the desired result.
\end{dow}

\begin{wn}
Assume that $f$ satisfies \rm{(A4$'$)}. Then there exists at most
one analytic solution  of \mbox{\rm(\ref{eq1.1})}.
\end{wn}
\begin{dow}
Follows from  Propositions \ref{uniqueness}, \ref{stw.12} and
\ref{stw4.6}.
\end{dow}

\begin{uw}
From Lemma \ref{lm.mr} it follows that the martingale $M$ from
Proposition \ref{stw.12} has the representation
\[
M_{t}=\int_{0}^{t} \nabla u(X_{\theta})\,d B_{\theta},
\quad t\ge 0,\quad P_{x}\mbox{-a.s.}
\]
for q.e. $x\in\Omega$.
\end{uw}

In the following three propositions by a solution of (\ref{eq1.1})
we mean either probabilistic solution or analytic solution.

Let us recall that by $T_r$ we denote the truncature operator
defined in the proof of Theorem \ref{main}.

\begin{stw}
\label{eq.w} Assume that $u$ is a solution of
\mbox{\rm(\ref{eq1.1})}  such that  $f(u)\in L^{1}(\Omega;m)$.
Then $T_{r}(u)\in H^{1}_{0}(\Omega)$ for every $r>0$, $u\in
W^{1,q}_{0}(\Omega)$ for $q\in[1,\frac{d}{d-1})$ and
\begin{equation}
\label{VI} \EE(u^{k},v)=(f^{k}(u),v)_{L^{2}(\Omega;m)}+\langle
v,\mu^{k}\rangle, \quad v\in C^{2}_{0}(\overline{\Omega}),\quad
k=1,\dots,n.
\end{equation}
\end{stw}
\begin{dow}
Since $f(u)\in L^{1}(\Omega;m)$ and $\mu\in\MM_{0,b}$,
$E_{x}\int_{0}^{\zeta}|f(u)(X_t)|\,dt<+\infty$  and
$E_{x}\int_{0}^{\zeta}d|A^{\mu}|_t<+\infty$ for q.e. $x\in\Omega$
(see the remark following (\ref{RD2})). Using this and the fact
that $u$ is of class (FD) we conclude from (\ref{eq.df.equiv})
that
\[
u(x)=E_{x}\int_{0}^{\zeta}f(u)(X_t)\,dt +E_x\int_{0}^{\zeta}
dA^{\mu}_t
\]
for q.e. $x\in\Omega$. From results proved in \cite{KR} it follows
now that $u$ is a solution of (\ref{eq1.1}) in the sense of
duality (see \cite{Stampacchia}). Consequently, $u$ satisfies
(\ref{VI}). Moreover, by \cite{Stampacchia}, $u\in
W^{1,q}_{0}(\Omega)$ with $1\le q<d/(d-1)$ while by \cite{DMOP},
$T_{r}(u)\in H^{1}_{0}(\Omega)$ for every $r>0$.
\end{dow}
\medskip

We do not know whether under (A4) the function $f(u)$ is
integrable. We can show, however, that it is integrable under
stronger that (A4) condition (A4$''$).

\begin{stw}
Let assumption  \mbox{\rm(A4$''$)} hold and let $u$ be a solution
of \mbox{\rm(\ref{eq1.1})}. Then
$f(u)\in L^{1}(\Omega;m)$, \mbox{\rm(\ref{VI})} is satisfied and
\[
\|f(u)\|_{L^{1}(\Omega;m)}\le \|\mu\|_{TV}.
\]
\end{stw}
\begin{dow}
Since $u$ is a solution of (\ref{eq1.1}), $u^{k}$ is a solution of
the scalar equation
\begin{equation*}
\left\{
\begin{array}{l} -\frac12\Delta u^{k}
=g^{k}(x,u^{k})+\mu^{k}\quad \mbox{in }\Omega,\medskip\\
u^{k}=0\quad \mbox{on }\partial\Omega,
\end{array}
\right.
\end{equation*}
where $g^{k}:\Omega\times\mathbb{R}\rightarrow\mathbb{R}$,
$g^{k}(x,t)=f^{k}(x,u^{1}(x),\dots,u^{k-1}(x),t,u^{k+1}(x),\dots,u^{n}(x))$.
Observe that by (A4$''$), $g^{k}(x,t)\cdot t\le 0$ for a.e.
$x\in\Omega$ and every $t\in\mathbb{R}$. Therefore by \cite[Lemma
2.3]{KR} (in \cite[Lemma 2.3]{KR} it is assumed that
$g^{k}(x,\cdot)$ is monotone but as a matter of fact in the proof
only the sign condition formulated above is used) and Proposition
\ref{stw.12}, for every $\tau\in\mathcal{T}$ such that
$0\le\tau<\zeta$ we have
\[
E_{x}\int_{0}^{\tau}|g^{k}(X_t,u^{k}(X_t))|\,dt \le
E_{x}|u^{k}(X_{\tau})| +E_{x}\int_{0}^{\tau}d|A^{\mu^{k}}|_t
\]
for q.e. $x\in\Omega$. Let $\{\tau_{k}\}\subset \mathcal{T}$ be
such that $0\le\tau_{k}<\zeta$ and $\tau_{k}\rightarrow \zeta$.
Replacing $\tau$ by $\tau_{k}$ in the above inequality, passing to
the limit and using the fact that $u$ is of class (FD) and
$u\in\mathcal{C}_{0}(\Omega)$ we get
\[
E_{x}\int_{0}^{\zeta}|g^{k}(X_t, u^{k}(X_t))|\,dt \le
E_{x}\int_{0}^{\zeta}d|A^{\mu^{k}}|_t.
\]
By the above and \cite[Lemma 5.4]{KR},
\[
\|f^{k}(u)\|_{L^{1}(\Omega;m)}=\|g^{k}(u_{k})\|_{L^{1}(\Omega;m)}\le
\|\mu^{k}\|_{TV},
\]
which implies the desired inequality. That (\ref{VI}) is satisfied
now follows from Proposition \ref{eq.w}.
\end{dow}
\medskip

Finally we show that if the right-hand side of (\ref{eq1.1})
satisfies the uniform angle condition (A5) then Stampacchia's
estimate (\ref{eq0.s}) holds for every solution of (\ref{eq1.1}).

\begin{stw}
Assume \mbox{\rm(A5)}. If $u$ is a solution of
\mbox{\rm(\ref{eq1.1})} then
\[
\|f(u)\|_{L^{1}(\Omega;m)}\le \alpha^{-1}\|\mu\|_{TV}.
\]
\end{stw}
\begin{dow}
Using the It\^o-Tanaka formula we get the first inequality in
(\ref{eq2.tf}). From this inequality, the fact that $u$ is of
class (FD), $u\in\mathcal{C}_{0}(\Omega)$ and (A5) one can
conclude that for every $\tau\in\mathcal{T}$ such that $0\le
\tau<\zeta$,
\[
E_{x}\int_{0}^{\tau}|f(u)(X_{t})|\,dt\le \alpha^{-1}
E_{x}\int_{0}^{\zeta}d|A^{\mu}|_{r}
\]
for q.e. $x\in\Omega$. The desired inequality now follows from
Fatou's lemma and \cite[Lemma 5.4]{KR}.
\end{dow}

\end{document}